\documentclass[10pt]{article}%
\usepackage{amsmath,amssymb,amscd,verbatim,amsthm}
\usepackage{amsmath}
\usepackage{amssymb}
\usepackage{amsfonts}
\usepackage{graphicx}%
\providecommand{\U}[1]{\protect\rule{.1in}{.1in}}
\begin{document}

\begin{center}
{\LARGE \textbf{New coincidence point and fixed point theorems for essential
distances and $e^{0}$-metrics}} \\[0.4in]

{\large \textbf{Wei-Shih Du}\footnote{{\small E-mail address:
wsdu@mail.nknu.edu.tw; Tel: +886-7-7172930 ext 6809; Fax: +886-7-6051061. }%
}\bigskip}

{\large Department of Mathematics, National Kaohsiung Normal
University,\newline Kaohsiung 82444, Taiwan }
\end{center}

\bigskip

\hrule\vspace{0.1cm}\bigskip

\noindent\textbf{Abstract: }In this paper, we establish some new fixed point
theorems and coincidence point theorems for essential distances and $e^{0}%
$-metrics which generalize and improve Berinde-Berinde's fixed point theorem,
Mizoguchi-Takahashi's fixed point theorem, Nadler's fixed point theorem and
Banach contraction principle and many known results in the literature.
\bigskip

\bigskip

\vspace{0.1cm}
\hrule\vspace{0.3cm} \bigskip

\noindent{\large \textbf{1. Introduction and preliminaries }}\bigskip

Let $(X,d)$ be a metric space. Denote by $\mathcal{N}(X)$ the class of all
nonempty subsets of $X$, $\mathcal{C}(X)$ the class of all nonempty closed
subsets of $X$ and $\mathcal{CB}(X)$ the family of all nonempty closed and
bounded subsets of $X$. A function $\mathcal{H}:\mathcal{CB}(X)\times
\mathcal{CB}(X)\rightarrow\lbrack0,\infty)$ defined by
\[
\mathcal{H}(A,B)=\text{max}\left\{  \sup_{x\in B}d(x,A)\text{,}\sup_{x\in
A}d(x,B)\right\}
\]
is said to be the Hausdorff metric on $\mathcal{CB}(X)$ induced by the metric
$d$ on $X$. A point $v$ in $X$ is a fixed point of a mapping $T$ if $v=Tv$
(when $T:X\rightarrow X$ is a single-valued mapping) or $v\in Tv$ (when
$T:X\rightarrow\mathcal{N}(X)$ is a multivalued mapping). The set of fixed
points of $T$ is denoted by $\mathcal{F}(T)$. Let $g:X\rightarrow X$ be a
self-mapping and $T:X\rightarrow\mathcal{N}(X)$ be a multivalued mapping. A
point $x$ in $X$ is said to be a \textit{coincidence point} of $g$ and $T$ if
$gx\in Tx$. The set of coincidence point of $g$ and $T$ is denoted by
$\mathcal{COP}(g,T)$. Throughout this paper, we denote by $%
%TCIMACRO{\U{2115} }%
%BeginExpansion
\mathbb{N}
%EndExpansion
$ and $%
%TCIMACRO{\U{211d} }%
%BeginExpansion
\mathbb{R}
%EndExpansion
$, the set of positive integers and real numbers, respectively.

In 2007, M. Berinde and V. Berinde [2] established the following interesting
fixed point theorem which generalizes Mizoguchi-Takahashi's fixed point
theorem [13], Nadler's fixed point theorem [14] and Banach contraction
principle [1].\bigskip

\noindent\textbf{Theorem 1.1 (M. Berinde and V. Berinde [2]).}\quad\textit{Let
}$(X,d)$\textit{ be a complete metric space, }$T:X\rightarrow\mathcal{CB}(X)$
\textit{be a multivalued mapping and }$L\geq0$\textit{. Assume that}
\[
\mathcal{H}(Tx,Ty)\leq\alpha(d(x,y))d(x,y)+Ld(y,Tx)\text{,}%
\]
\textit{for all }$x,y\in X$\textit{, where }$\alpha$\textit{ is an}
$\mathcal{MT}$\textit{-function, that is, }$\alpha$\textit{ is a function from
}$[0,\infty)$\textit{ into }$[0,1)$\textit{ satisfying }$\limsup_{s\rightarrow
t+0}\alpha(s)<1$\textit{ for all }$t\in\lbrack0,\infty)$\textit{. Then there
exists }$v\in X$\textit{ such that }$v\in Tv$\textit{.}\bigskip

Let $(X,d)$ be a metric space. Recall that a function $\kappa:X\times
X\rightarrow\lbrack0,\infty)$ is said to be a $\tau$\textit{-function} [3, 4,
6, 7, 10, 12], introduced and studied by Lin and Du, if the following
conditions hold:

\begin{enumerate}
\item[($\tau1$)] $\kappa(x,z)\leq\kappa(x,y)+\kappa(y,z)$ for all $x,y,z\in X$;

\item[($\tau2$)] if $x\in X$ and $\{y_{n}\}$ in $X$ with $\lim_{n\rightarrow
\infty}y_{n}=y$ such that $\kappa(x,y_{n})\leq M$ for some $M=M(x)>0$, then
$\kappa(x,y)\leq M$;

\item[($\tau3$)] for any sequence $\{x_{n}\}$ in $X$ with $\lim_{n\rightarrow
\infty}\sup\{\kappa(x_{n},x_{m}):m>n\}=0$, if there exists a sequence
$\{y_{n}\}$ in $X$ such that $\lim_{n\rightarrow\infty}\kappa(x_{n},y_{n})=0$,
then $\lim_{n\rightarrow\infty}\rho(x_{n},y_{n})=0$;

\item[($\tau4$)] for $x,$ $y,$ $z\in X$, $\kappa(x,y)=0$ and $\kappa(x,z)=0$
imply $y=z$.
\end{enumerate}

It is obvious that the metric $d$ is a $w$-distance [3, 4, 6, 7, 9-12] and any
$w$-distance is a $\tau$-function, but the converse is not true; see [3, 10,
12] for more detail.

The following result is very crucial in our proofs.\bigskip

\noindent\textbf{Lemma 1.1 (see [8, Lemma 1.1]).}\quad\textit{If condition
}$(\tau4)$\textit{\ is weakened to the following condition }$(\tau4)^{\prime
}:$

\begin{enumerate}
\item[$(\tau4)^{\prime}$] \textit{for any }$x\in X$\textit{\ with }%
$\kappa(x,x)=0$\textit{, if }$\kappa(x,y)=0$\textit{\ and }$\kappa
(x,z)=0$\textit{, then }$y=z$\textit{,}
\end{enumerate}

\noindent\textit{then }$(\tau3)$\textit{\ implies }$(\tau4)^{\prime}%
$\textit{.}\bigskip

In 2016, Du [8] introduced the concept of essential distance as
follows.\bigskip

\noindent\textbf{Definition 1.1 (see [8, Definition 1.2]).}\quad Let $(X,d)$
be a metric space. A function $\kappa:X\times X\rightarrow\lbrack0,+\infty)$
is called an \textit{essential distance} (abbreviated as "$e$%
-\textit{distance}") if conditions {$(\tau1)$, $(\tau2)$ and $(\tau3)$
hold.}\bigskip

\noindent\textbf{Remark 1.1.}\quad Clearly any $\tau$-function is an
$e$-distance. By Lemma 1.1, we know that if $\kappa$ is an $e$-distance, then
condition $(\tau4)^{\prime}$ holds.\bigskip

The following known result is very crucial in this paper.\bigskip

\noindent\noindent\textbf{Lemma 1.2 (see [4, Lemma 2.1]).}\quad\textit{Let
}$(X,d)$\textit{\ be a metric space and }$\kappa:X\times X\rightarrow
\lbrack0,+\infty)$\textit{\ be a function. Assume that }$\kappa$%
\textit{\ satisfies the condition }$(\tau3)$\textit{. If a sequence }%
$\{x_{n}\}$\textit{\ in }$X$\textit{\ with }$\lim_{n\rightarrow\infty}%
\sup\{\kappa(x_{n},x_{m}):m>n\}=0$\textit{, then }$\{x_{n}\}$\textit{\ is a
Cauchy sequence in }$X$\textit{.}\bigskip

In 2016, Du introduced the concept of $\mathcal{MT}(\lambda)$-function [7] as
follows.\bigskip

\noindent\textbf{Definition 1.2.}\quad Let $\lambda>0$. A function $\mu:$
$[0,+\infty)\rightarrow$ $[0,\lambda)$ is said to be an $\mathcal{MT}%
(\lambda)$-\textit{function} [7] if $\limsup\limits_{s\rightarrow t^{+}}%
\mu(s)<\lambda$ for all $t\in\lbrack0,\infty)$. As usual, we simply write
"$\mathcal{MT}$-function (see [3-6])" instead of "$\mathcal{MT}(1)$%
-function".\bigskip

In [7], Du established the following useful characterizations of
$\mathcal{MT}(\lambda)$-functions.\bigskip

\noindent\textbf{Theorem 1.2 (see [7, Theorem 2.4]).}\quad\textit{Let
}$\lambda>0$\textit{\ and let }$\mu:$\textit{\ }$[0,+\infty)\rightarrow
$\textit{\ }$[0,\lambda)$\textit{\ be a function. Then the following
statements are equivalent.}

\begin{enumerate}
\item[(1)] $\mu$\textit{\ is an }$\mathcal{MT}(\lambda)$\textit{-function.}

\item[(2)] $\lambda^{-1}\mu$\textit{\ is an }$\mathcal{MT}$\textit{-function.}

\item[(3)] \textit{For each }$t\in\lbrack0,\infty)$\textit{, there exist }%
$\xi_{t}^{(1)}\in\lbrack0,\lambda)$\textit{\ and }$\epsilon_{t}^{(1)}%
>0$\textit{\ such that }$\mu(s)\leq\xi_{t}^{(1)}$\textit{\ for all }%
$s\in(t,t+\epsilon_{t}^{(1)})$\textit{.}

\item[(4)] \textit{For each }$t\in\lbrack0,\infty)$\textit{, there exist }%
$\xi_{t}^{(2)}\in\lbrack0,\lambda)$\textit{\ and }$\epsilon_{t}^{(2)}%
>0$\textit{\ such that }$\mu(s)\leq\xi_{t}^{(2)}$\textit{\ for all }%
$s\in\lbrack t,t+\epsilon_{t}^{(2)}]$\textit{.}

\item[(5)] \textit{For each }$t\in\lbrack0,\infty)$\textit{, there exist }%
$\xi_{t}^{(3)}\in\lbrack0,\lambda)$\textit{\ and }$\epsilon_{t}^{(3)}%
>0$\textit{\ such that }$\mu(s)\leq\xi_{t}^{(3)}$\textit{\ for all }%
$s\in(t,t+\epsilon_{t}^{(3)}]$\textit{.}

\item[(6)] \textit{For each }$t\in\lbrack0,\infty)$\textit{, there exist }%
$\xi_{t}^{(4)}\in\lbrack0,\lambda)$\textit{\ and }$\epsilon_{t}^{(4)}%
>0$\textit{\ such that }$\mu(s)\leq\xi_{t}^{(4)}$\textit{\ for all }%
$s\in\lbrack t,t+\epsilon_{t}^{(4)})$\textit{.}

\item[(7)] \textit{For any nonincreasing sequence }$\{x_{n}\}_{n\in%
%TCIMACRO{\U{2115} }%
%BeginExpansion
\mathbb{N}
%EndExpansion
}$\textit{\ in }$[0,\infty)$\textit{, we have }$0\leq\sup\limits_{n\in%
%TCIMACRO{\U{2115} }%
%BeginExpansion
\mathbb{N}
%EndExpansion
}\mu(x_{n})<\lambda$\textit{.}

\item[(8)] \textit{For any strictly decreasing sequence }$\{x_{n}\}_{n\in%
%TCIMACRO{\U{2115} }%
%BeginExpansion
\mathbb{N}
%EndExpansion
}$\textit{\ in }$[0,\infty)$\textit{, we have }$0\leq\sup\limits_{n\in%
%TCIMACRO{\U{2115} }%
%BeginExpansion
\mathbb{N}
%EndExpansion
}\mu(x_{n})<\lambda$\textit{.}

\item[(9)] \textit{For any eventually nonincreasing sequence }$\{x_{n}\}_{n\in%
%TCIMACRO{\U{2115} }%
%BeginExpansion
\mathbb{N}
%EndExpansion
}$\textit{\ }(\textit{i.e. there exists }$\ell\in%
%TCIMACRO{\U{2115} }%
%BeginExpansion
\mathbb{N}
%EndExpansion
$\textit{\ such that }$x_{n+1}\leq x_{n}$\textit{\ for all }$n\in%
%TCIMACRO{\U{2115} }%
%BeginExpansion
\mathbb{N}
%EndExpansion
$\textit{\ with }$n\geq\ell$) \textit{in }$[0,\infty)$\textit{, we have
}$0\leq\sup\limits_{n\in%
%TCIMACRO{\U{2115} }%
%BeginExpansion
\mathbb{N}
%EndExpansion
}\mu(x_{n})<\lambda$\textit{.}

\item[(10)] \textit{For any eventually strictly decreasing sequence }%
$\{x_{n}\}_{n\in%
%TCIMACRO{\U{2115} }%
%BeginExpansion
\mathbb{N}
%EndExpansion
}$\textit{\ }(\textit{i.e. there exists }$\ell\in%
%TCIMACRO{\U{2115} }%
%BeginExpansion
\mathbb{N}
%EndExpansion
$\textit{\ such that }$x_{n+1}<x_{n}$\textit{\ for all }$n\in%
%TCIMACRO{\U{2115} }%
%BeginExpansion
\mathbb{N}
%EndExpansion
$\textit{\ with }$n\geq\ell$) \textit{in }$[0,\infty)$\textit{, we have
}$0\leq\sup\limits_{n\in%
%TCIMACRO{\U{2115} }%
%BeginExpansion
\mathbb{N}
%EndExpansion
}\mu(x_{n})<\lambda$\textit{.}\bigskip
\end{enumerate}

Let $\kappa$ be an $e$-function. For each $x\in X$ and $C\subseteq X$, denote
by
\[
\kappa(x,C)=\inf_{y\in C}\kappa(x,y).
\]
\medskip

\noindent\textbf{Lemma 1.3 (see [3, Lemma 1.2].}\quad\textit{Let }%
$C$\textit{\ be a closed subset of a metric space }$(X,d)$\textit{\ and
}$\kappa$\textit{\ be a function satisfying the condition }$(\tau3)$\textit{.
Suppose that there exists }$z\in X$\textit{\ such that }$\kappa(z,z)=0$%
\textit{. Then }$\kappa(z,C)=0$\textit{\ if and only if }$z\in C$\textit{.}
\bigskip

Now, we introduce the concepts of $e^{0}$-distance and $e^{0}$-metric.\bigskip

\noindent\textbf{Definition 1.3.}\quad Let $(X,d)$ be a metric space. A
function $\kappa:X\times X\rightarrow\lbrack0,\infty)$ is called an $e^{0}%
$-\textit{distance} if it is an $e$-distance on $X$ with $\kappa(x,x)=0$ for
all $x\in X$.\bigskip

\noindent\noindent\textbf{Definition 1.4.}\quad Let $(X,d)$ be a metric space
and $\kappa$ be an $e^{0}$-distance. For any $A$, $B\in\mathcal{CB}(X)$,
define a function\ $\mathcal{D}_{\kappa}:\mathcal{CB}(X)\times\mathcal{CB}%
(X)\rightarrow\lbrack0,+\infty)$ by
\[
\mathcal{D}_{\kappa}(A,B)=\text{max}\{\xi_{\kappa}(A,B),\xi_{\kappa
}(B,A)\}\text{,}%
\]
where $\xi_{\kappa}(A,B)=\sup_{x\in A}\kappa(x,B)$, then $\mathcal{D}_{\kappa
}$ is said to be the $e^{0}$-$metric$ on $\mathcal{CB}(X)$ induced by $\kappa
$.\bigskip

Clearly, any Hausdorff metric is an $e^{0}$-metric, but the reverse is not
true. It is not difficult to prove the following theorem. \bigskip

\noindent\textbf{Theorem 1.3.}\quad\textit{Let }$(X,d)$\textit{\ be a metric
space and }$\mathcal{D}_{\kappa}$\textit{\ be an }$e^{0}$\textit{-metric
defined as in Def. 1.4 on }$\mathcal{CB}(X)$\textit{\ induced by an }$e^{0}%
$\textit{-distance }$\kappa$\textit{. Then for }$A$\textit{, }$B$\textit{,
}$C\in\mathcal{CB}(X)$\textit{, the following hold:}

\begin{enumerate}
\item[(i)] $\xi_{\kappa}(A,B)=0$\textit{\ }$\Longleftrightarrow$%
\textit{\ }$A\subseteq B$\textit{;}

\item[(ii)] $\xi_{\kappa}(A,B)\leq\xi_{\kappa}(A,C)+\xi_{\kappa}%
(C,B)$\textit{;}

\item[(iii)] \textit{Every }$e^{0}$\textit{-metric }$\mathcal{D}_{\kappa}%
$\textit{\ is a metric on }$\mathcal{CB}(X)$\textit{.}
\end{enumerate}

\bigskip

In this paper, some new fixed point theorems and coincidence point theorems
for essential distances and $e^{0}$-metrics are established. Our new results
generalize and improve Berinde-Berinde's fixed point theorem,
Mizoguchi-Takahashi's fixed point theorem, Nadler's fixed point theorem,
Banach contraction principle and many known results in the literature.
\bigskip\bigskip\bigskip

\noindent{\large \textbf{2. New fixed point and coincidence point theorems}%
}\bigskip

In this section, we first establish the following new fixed point theorem for
$e^{0}$-distances which extends and generalizes Berinde-Berinde's fixed point
theorem, Mizoguchi-Takahashi's fixed point theorem, Nadler's fixed point
theorem and Banach contraction principle.\bigskip

\noindent\textbf{Theorem 2.1.}\quad\textit{Let }$(X,d)$\textit{ be a complete
metric space, }$\kappa$\textit{ be a }$e^{0}$\textit{-distance and
}$T:X\rightarrow\mathcal{C}(X)$\textit{ be a multivalued mapping. Suppose
that}

\begin{enumerate}
\item[$(S1)$] \textit{there exists an} $\mathcal{MT}$\textit{-function }$\mu
:$\textit{ }$[0,\infty)\rightarrow$\textit{ }$[0,1)$\textit{ such that for
each }$x\in X$\textit{, if }$y\in Tx$\textit{ with }$y\neq x$\textit{ then
there exists }$z\in Ty$\textit{ such that }%
\[
\kappa(y,z)\leq\mu(k(x,y))\kappa(x,y)\text{;}%
\]

\item[$(S2)$] $T$\textit{ further satisfies one of the following conditions:}

\begin{enumerate}
\item[(H1)] $T$\textit{ is closed;}

\item[(H2)] \textit{the function }$f:X\rightarrow\lbrack0,\infty)$\textit{
defined by }$f(x)=\kappa(x,Tx)$\textit{ is l.s.c.;}

\item[(H3)] \textit{the function }$g:X\rightarrow\lbrack0,\infty)$\textit{
defined by }$g(x)=d(x,Tx)$\textit{ is l.s.c.;}

\item[(H4)] \textit{for any sequence }$\{x_{n}\}$\textit{ in }$X$\textit{ with
}$x_{n+1}\in Tx_{n}$\textit{, }$n\in%
%TCIMACRO{\U{2115} }%
%BeginExpansion
\mathbb{N}
%EndExpansion
$\textit{ and }$\lim_{n\rightarrow\infty}x_{n}=v$\textit{, we have }%
$\lim_{n\rightarrow\infty}\kappa(x_{n},Tv)=0$\textit{;}

\item[(H5)] $\inf\{\kappa(x,z)+\kappa(x,Tx):$\textit{ }$x\in X\}>0$\textit{
for every }$z\notin\mathcal{F}(T)$\textit{.}
\end{enumerate}
\end{enumerate}

\noindent\textit{Then }$\mathcal{F}(T)\neq\emptyset$\textit{.} \medskip

\noindent\textbf{Proof.}\quad Since $\kappa$ is a $e^{0}$-distance, by Lemma
1.1, we know that for $x,y\in X$, $\kappa(x,y)=0$\textit{ }%
$\Longleftrightarrow$\textit{ }$x=y$. Following a similar argument as the
proof of [6, Lemma 3.1], we can prove the conclusion.\hspace{\fill} $\Box
$\bigskip

The following result is an immediate consequence of Theorem 2.1.\bigskip

\noindent\noindent\textbf{Theorem 2.2.}\quad\textit{Let }$(X,d)$\textit{ be a
complete metric space, }$\kappa$\textit{ be a }$e^{0}$\textit{-distance,
}$T:X\rightarrow C(X)$\textit{ be a multivalued mapping and }$\mu:$\textit{
}$[0,\infty)\rightarrow$\textit{ }$[0,1)$\textit{ be an }$\mathcal{MT}%
$\textit{-function. Suppose that the condition }$(S2)$\textit{ as in Theorem
2.1 holds and further assume that}

\begin{enumerate}
\item[$(S3)$] \textit{for each }$x\in X$\textit{, }$\kappa(y,Ty)\leq\mu
(\kappa(x,y))\kappa(x,y)$\textit{ for all }$y\in Tx$\textit{.}
\end{enumerate}

\noindent\textit{Then }$\mathcal{F}(T)\neq\emptyset$\textit{.\bigskip}

By applying Theorem 2.2, we establishe the following existence theorem of
coincidence point and fixed point.\bigskip

\noindent\textbf{Theorem 2.3.}\quad\textit{Let }$(X,d)$\textit{ be a complete
metric space, }$\kappa$\textit{ be a }$e^{0}$\textit{-distance, }%
$T:X\rightarrow\mathcal{C}(X)$ \textit{be a multivalued mapping, }%
$\varphi:X\rightarrow X$\textit{ be a continuous self-mapping, }$\mu:$\textit{
}$[0,\infty)\rightarrow$\textit{ }$[0,1)$\textit{ be an} $\mathcal{MT}%
$\textit{-function and }$L\geq0$\textit{. Suppose that the condition }%
$(S2)$\textit{ as in Theorem 2.1 holds and further assume}

\begin{enumerate}
\item[$(S4)$] $Tx$\textit{ is }$\varphi$\textit{-invariant (i.e. }%
$\varphi(Tx)\subseteq Tx$\textit{) for each }$x\in X$\textit{;}

\item[$(S5)$] $\kappa(y,Ty)\leq\mu(\kappa(x,y))\kappa(x,y)+L\kappa(\varphi
y,Tx)$\textit{ for all }$x,y\in X$\textit{.}
\end{enumerate}

\noindent\textit{Then }$\mathcal{COP}(\varphi,T)\cap\mathcal{F}(T)\neq
\emptyset$\textit{.\bigskip}

The following existence theorem of coincidence point and fixed point
for\ $e^{0}$-distances and $e^{0}$-metrics is an immediate consequence of
Theorem 2.3.\bigskip

\noindent\textbf{Theorem 2.4.}\quad\textit{Let }$(X,d)$\textit{ be a complete
metric space, }$\kappa$\textit{ be a }$e^{0}$\textit{-distance,}
$\mathcal{D}_{\kappa}$ \textit{be a }$e^{0}$\textit{-metric on} $\mathcal{CB}%
(X)$, $T:X\rightarrow\mathcal{CB}(X)$ \textit{be a multivalued mapping,
}$\varphi:X\rightarrow X$\textit{ be a continuous self-mapping, }$\mu
:$\textit{ }$[0,\infty)\rightarrow$\textit{ }$[0,1)$\textit{ be an}
$\mathcal{MT}$\textit{-function and }$L\geq0$\textit{. Suppose that the
conditions }$(S2)$\textit{ and }$(S4)$\textit{ hold and further assume}

\begin{enumerate}
\item[$(S6)$] $\mathcal{D}_{\kappa}(Tx,Ty)\leq\mu(\kappa(x,y))\kappa
(x,y)+L\kappa(\varphi y,Tx)$\textit{ for all }$x,y\in X$\textit{.}
\end{enumerate}

\noindent\textit{Then }$\mathcal{COP}(\varphi,T)\cap\mathcal{F}(T)\neq
\emptyset$\textit{.}\bigskip

\noindent\textbf{Remark 2.1.}\quad Theorems 2.1-2.4 all improve and generalize
Berinde-Berinde's fixed point theorem, Mizoguchi-Takahashi's fixed point
theorem, Nadler's fixed point theorem, Banach contraction principle and some
results in [3, 5, 6] and references therein.\bigskip\bigskip\bigskip

\noindent{\large \textbf{Acknowledgments:}} The author is supported by Grant
No. MOST 107-2115-M-017-004-MY2 of the Ministry of Science and Technology of
the Republic of China.\bigskip\bigskip

\noindent{\large \textbf{References}}

\begin{enumerate}
\item[{[1]}] S. Banach, \textit{Sur les op\'{e}rations dans les ensembles
abstraits et leurs applications aux \'{e}quations int\'{e}grales}, Fundamenta
Mathematicae 3 (1922) 133-181.

\item[{[2]}] M. Berinde, V. Berinde, \textit{On a general class of
multi-valued weakly Picard mappings}, J. Math. Anal. Appl. 326 (2007) 772-782.

\item[{[3]}] W.-S. Du, \textit{Some new results and generalizations in metric
fixed point theory}, Nonlinear Anal. 73 (2010) 1439-1446.

\item[{[4]}] W.-S. Du, \textit{Critical point theorems for nonlinear dynamical
systems and their applications}, Fixed Point Theory and Applications (2010),
Article ID 246382, doi:10.1155/2010/246382.

\item[{[5]}] W.-S. Du, \textit{On coincidence point and fixed point theorems
for nonlinear multivalued maps}. Topol. Appl. 159 (2012) 49--56.

\item[{[6]}] W.-S. Du, S.-X. Zheng, \textit{Nonlinear conditions for
coincidence point and fixed point theorems}, Taiwan. J. Math. 16(3) (2012) 857--868.

\item[{[7]}] W.-S. Du, \textit{New existence results of best proximity points
and fixed points for }$MT(\lambda)$\textit{-functions with applications to
differential equations}, Linear Nonlinear Anal. 2(2) (2016) 199--213.

\item[{[8]}] W.-S. Du, \textit{On generalized Caristi's fixed point theorem
and its equivalence}, Nonlinear Anal. Differ. Equ. 4(13) (2016) 635--644.

\item[{[9]}] O. Kada, T. Suzuki, W. Takahashi, \textit{Nonconvex minimization
theorems and fixed point theorems in complete metric spaces}, Math. Japon. 44
(1996) 381-391.

\item[{[10]}] L.-J. Lin, W.-S. Du, \textit{Ekeland's variational principle,
minimax theorems and existence of nonconvex equilibria in complete metric
spaces}, J. Math. Anal. Appl. 323 (2006) 360-370.

\item[{[11]}] L.-J. Lin, W.-S. Du, \textit{Some equivalent formulations of
generalized Ekeland's variational principle and their applications}, Nonlinear
Anal. 67 (2007) 187-199.

\item[{[12]}] L.-J. Lin, W.-S. Du, \textit{On maximal element theorems,
variants of Ekeland's variational principle and their applications}, Nonlinear
Anal. 68 (2008) 1246-1262.

\item[{[13]}] N. Mizoguchi, W. Takahashi, \textit{Fixed point theorems for
multivalued mappings on complete metric spaces}, J. Math. Anal. Appl. 141
(1989) 177-188.

\item[{[14]}] S.B. Nadler, Jr., \textit{Multi-valued contraction mappings},
Pacific J. Math. 30 (1969) 475-488.

%%%%%%%%%%%%%%%%%%%%%%%%%%%%%%%%11.25
\newpage
\end{enumerate}

\end{document}